\documentclass[letterpaper, 10 pt, conference]{ieeeconf}

\usepackage{amssymb,amsmath,color,amsfonts}
\usepackage{graphicx}
\usepackage{caption}
\usepackage{cite}
\IEEEoverridecommandlockouts                              
\newcommand*{\QEDB}{\hfill\ensuremath{\square}}%

\newtheorem{assumption}{Assumption}

\newtheorem{lemma}{Lemma}
\newtheorem{theorem}{Theorem}

\newtheorem{remark}{Remark}

\begin{document}

\title{\bf Fixed-Time Nash Equilibrium Seeking in Non-Cooperative Games}

\author{Jorge I. Poveda\thanks{ Jorge I. Poveda is with the Department of Electrical, Computer, and Energy Engineering, University of Colorado, Boulder, CO, 80309. Email:\{jorge.poveda@colorado.edu.\}}, Miroslav Krsti\'{c}\thanks{Miroslav Krsti\'{c} is with the Department of Mechanical and Aerospace Engineering, University of California, La Jolla, CA}, Tamer Ba\c{s}ar \thanks{Tamer Ba\c{s}ar is with the Department of Electrical and Computer Engineering, University of Illinois at Urbana-Champaign, IL.}  \thanks{Research supported in part by AFOSR grants FA9550-18-1-0246 and FA9550-19-1-0353, NSF grants CNS-1947613 and ECCS-1508757, and a CU Boulder ASIRT grant. Submitted on March 30, 2020.}}

\maketitle
\begin{abstract}
We introduce a novel class of Nash equilibrium seeking dynamics for non-cooperative games with a finite number of players, where the convergence to the Nash equilibrium is bounded by a $\mathcal{K}\mathcal{L}$ function with a settling time that can be upper bounded by a positive constant that is independent of the initial conditions of the players, and which can be prescribed a priori by the system designer. The dynamics are model-free, in the sense that the mathematical forms of the cost functions of the players are unknown. Instead, in order to update its own action, each player needs to have access only to real-time evaluations of its own cost, as well as to auxiliary states of neighboring players characterized by a communication graph. Stability and convergence properties are established for both potential games and strongly monotone games. Numerical examples are presented to illustrate our theoretical results.
\end{abstract}
\section{Introduction}
In many engineering and socio-technical systems, there exists inherent competition between different entities or subsystems who aim to maximize their individual payoffs by controlling their own actions. Since, in general, the payoffs of the subsystems also depend on the actions of each other, this setting describes a standard non-cooperative game, where the notion of \emph{Nash equilibrium} has played an important role during the last decades by providing a rigorous mathematical characterization of operational points where players have no incentive to deviate \cite{BasarDNGT}. 

In the controls and optimization literatures, several iterative algorithms have been developed to compute the Nash equilibrium of a given game under different assumptions on the information available to the players \cite{BasarAsynchronous}. In \cite{Frihauf12a}, it was shown that model-free dynamics inspired by extremum seeking ideas can also be used to achieve Nash equilibrium seeking (NES) in non-cooperative games with static and dynamic players. Similar ideas have been pursued in \cite{Kutadinata_journal}, \cite{PoQu15}, \cite{Hu15}, \cite{PoTe16} and \cite{Stipanovik_ES_TAC}, under different assumptions on the structure of the game and the communication requirements of the algorithm. However, while several algorithms have been able to solve the Nash equilibrium seeking problem in a model-free way, achieving desirable rates of convergence has been a persistent challenge in the design of model-free Nash seeking algorithms. Indeed, existing NES dynamics based on extremum seeking have been able to achieve only asymptotic results, where the bound that dominates the convergence of the actions of the players is either an exponential function \cite{Frihauf12a} or just a general $\mathcal{K}\mathcal{L}$ bound \cite{Kutadinata_journal}, \cite{PoQu15}. In these cases, and unless the gain of the dynamics is continuously increased, the convergence time of the algorithms will grow unbounded as the compact set of initial conditions also grows.

On the other hand, during recent years there have been several new results in the context of \emph{fixed-time} optimization and \emph{fixed-time} control \cite{Fixed_timeTAC,fixedtimebook,fixed_time}. This class of algorithms can achieve finite-time convergence to a desired target, with a finite-time that can be upper bounded by a constant that is independent of the initial conditions.  However, most of the results developed so far in the context of fixed-time multi-agent learning and optimization are applicable only to systems for which a precise mathematical model is available \cite{fixed_time}, or to single-agent learning problems with no spatial information constraints \cite{JournalPovedaKrstic}. Motivated by this background, in this paper, we present the first model-free NES dynamics with the (semi-global practical) \emph{fixed-time} convergence property for non-cooperative games with a finite number of players. To be more precise, we propose a new \emph{model-free} learning algorithm which guarantees that the actions of the players converge to a neighborhood of the Nash equilibrium of the game, with a convergence bound characterized by a $\mathcal{K}\mathcal{L}$ function having a uniformly bounded settling time. Moreover, the bound on the settling time can be prescribed a priori under a mild knowledge of the monotonicity properties of the game. Given that these types of dynamics are necessarily non-Lipschitz continuous, they cannot be studied by using standard tools of smooth extremum seeking control, e.g., \cite{krstic_book_ESC}. Instead, we use averaging tools for non-smooth extremum seeking dynamics, see e.g., \cite{WangAveraging,zero_order_poveda_Lina}. Numerical simulations are presented to illustrate the results.

The rest of this paper is organized as follows. Section \ref{sec_notation} presents the  preliminaries. Section \ref{sec_problem} presents the main NES dynamics and the main results. The proofs are presented in Section \ref{sec_proofs}. Numerical results are presented in Section \ref{sec_numerical}, and Section \ref{section:conclusions} ends with the conclusions.
%
%
%
\section{Notation and Definitions}
\label{sec_notation}
Given a compact set $\mathcal{A}\subset\mathbb{R}^n$ and a vector $z\in\mathbb{R}^n$, we use $|z|_{\mathcal{A}}:=\min_{s\in\mathcal{A}}\|z-s\|_2$ to denote the minimum distance of $z$ to $\mathcal{A}$.  We use $\mathbb{S}^1:=\{z\in\mathbb{R}^2:z^2_1+z_2^2=1\}$ to denote the unit circle in $\mathbb{R}^2$, and $\mathbb{T}^N:=\mathbb{S}^1\times \mathbb{S}^1\times\ldots\times \mathbb{S}^1$ to denote the Cartesian product of order $N$ of the set $\mathbb{S}_1$. We also use $r\mathbb{B}$ to denote a closed ball in the Euclidean space, of radius $r>0$, and centered at the origin. We use $I_n\in\mathbb{R}^{n\times n}$ to denote the identity matrix. A function $\beta:\mathbb{R}_{\geq0}\times\mathbb{R}_{\geq0}\to\mathbb{R}_{\geq0}$ is of class $\mathcal{K}\mathcal{L}$ if it is nondecreasing in its first argument, nonincreasing in its second argument, $\lim_{r\to0^+}\beta(r,s)=0$ for each $s\in\mathbb{R}_{\geq0}$, and  $\lim_{s\to\infty}\beta(r,s)=0$ for each $r\in\mathbb{R}_{\geq0}$. A function $\tilde{\beta}$ is of class $\mathcal{K}\mathcal{L}_{\mathcal{T}}$ \cite{RiosTeel} if $\tilde{\beta}\in\mathcal{K}\mathcal{L}$, and additionally, for each $r>0$  there exists a function $T:\mathbb{R}\to\mathbb{R}_{\geq0}$, called the settling time, such that $\tilde{\beta}(r,s)=0$ for all $s>T(r)$. When there exists $m>0$ such that $T(r)<m$ for all $r>0$, we say that $\tilde{\beta}$ has the \emph{fixed-time convergence} property. Given a vector-valued function $H:\mathbb{R}^N\to\mathbb{R}^N$, we define its \emph{pseudo-gradient} $G$ as the vector-valued function  $G(x):=[\frac{\partial H_1}{\partial x_1},\frac{\partial H_2}{\partial x_2},\ldots,\frac{\partial H_N}{\partial x_N}]^\top$. For the purpose of algorithmic analysis, in this paper we will work with constrained dynamical systems. In particular, we will consider ordinary differential equations (ODEs) with state $x\in\mathbb{R}^n$, and dynamics
\begin{equation}\label{ODE}
x\in C,~~~~\dot{x}=F(x),
\end{equation}
where $F:\mathbb{R}^n\to\mathbb{R}^n$ is a continuous function, and $C\subset\mathbb{R}^n$ is a closed set. A solution to system \eqref{ODE} is a continuously differentiable function $x:\text{dom}(x)\to\mathbb{R}^n$ that satisfies: a) $x(0)\in C$; b) $x(t)\in C$ for all $t\in\text{dom}(x)$; and c) $\dot{x}(t)=F(x(t))$ for all $t\in\text{dom}(x)$. A solution is said to be complete if $\text{dom}(x)=[0,\infty)$. Given a compact set $\mathcal{A}\subset C$, system \eqref{ODE} is said to render $\mathcal{A}$ uniformly globally asymptotically stable (UGAS) if there exists a class $\mathcal{K}\mathcal{L}$ function $\beta$ such that every solution of \eqref{ODE} satisfies  $|x(t)|_{\mathcal{A}}\leq \beta(|x(0)|_{\mathcal{A}},t$), for all $t\in\text{dom}(x)$. In this paper, we will also consider $\varepsilon$-perturbed or $\varepsilon$-parameterized dynamical systems of the form
\begin{equation}\label{perturbed_ode}
x\in C,~~~\dot{x}=F_{\varepsilon}(x),
\end{equation}
where $F_{\varepsilon}$ is a continuous function parameterized by a positive constant number $\varepsilon>0$. For these systems, we will study semi-global practical stability properties. In particular, a compact set $\mathcal{A}\subset C$ is said to be $\beta$-Semi-Globally Practically Asymptotically Stable ($\beta$-SGPAS) as $\varepsilon\to0^+$, if there exists a class $\mathcal{K}\mathcal{L}$ function $\beta$ such that for each pair $\delta>\nu>0$ there exists $\varepsilon^*>0$ such that for all $\varepsilon\in(0,\varepsilon^*)$ every solution of \eqref{perturbed_ode} with $|x(0)|_{\mathcal{A}}\leq \delta$ satisfies 
\begin{equation}
|x(t)|_{\mathcal{A}}\leq \beta(|x(0)|_{\mathcal{A}},t)+\nu,~~~\forall~t\in\text{dom}(x).
\end{equation}
The notion of SGPAS can be extended to systems that depend on multiple parameters $\varepsilon=[\varepsilon_1,\varepsilon_2,\ldots,\varepsilon_{\ell}]^\top$. In this case, and with some abuse of notation, we say that the system \eqref{perturbed_ode} renders the set $\mathcal{A}$ $\beta$-SGPAS as $(\varepsilon_{\ell},\ldots,\varepsilon_2,\varepsilon_{1})\to0^+$, where the parameters are tuned in order starting from $\varepsilon_1$. 

\section{Fixed-Time Nash Equilibrium Seeking Dynamics}
\label{sec_problem}
We consider a non-cooperative game with $N$ players, where each player can control only its own action $u_i\in\mathbb{R}$. For each player $i\in\mathcal{V}:=\{1,2,3,\ldots,N\}$, we denote by $u_{-i}\in\mathbb{R}^{N-1}$ the actions of the other players, and $J_i(u_i,u_{-i})$ as the cost function of the $i^{th}$ player, which in this paper is assumed to be a real-valued function, i.e., $J_i:\mathbb{R}^N\to\mathbb{R}$. The goal of the players is to converge to a Nash equilibrium of the game, which is an action profile $u^*\in\mathbb{R}^N$ that satisfies the condition  $J_i(u_i^*,u^*_{-i})= \inf_{u_i\in\mathbb{R}}J_i(u_i,u_{-i}^*)$, for all $i\in\mathcal{V}$. In this paper, we assume that the cost functions $\{J_i\}_{i=1}^N$ are such that the underlying game admits a unique Nash equilibrium\footnote{For the kinds of games that are our focus in this paper (potential games and strongly monotone games)  the condition of existence of a unique Nash equilibrium is satisfied, as discussed later.}, and to reach that equilibrium in a distributed way we are interested in strategies that have two main properties: (a) Their convergence can be characterized by class $\mathcal{K}\mathcal{L}_{\mathcal{T}}$ functions $\tilde{\beta}$ that have the fixed-time convergence property; (b) We assume that the mathematical forms of the cost functions $J_i$ and their partial derivatives are unknown, and each player has access only to local real-time evaluations of its own cost, as well as to auxiliary states of neighboring players characterized by a communication graph $\mathcal{G}$.
%
%

In order to achieve Nash equilibrium seeking under these two requirements, we will consider a class of model-free fixed-time Nash equilibrium seeking (FxTNES) dynamics based on fixed-time extremum seeking ideas \cite{PovedaKrsticACC20}, \cite{PovedaKrsticIFACWC20}.  We will focus on non-cooperative games that can be modeled either as potential games or strongly monotone games. To simplify our presentation, we will make the following assumption on the communication network of the game.
\vspace{0.1cm}

\begin{assumption}\label{graph_assumption}
The communication graph $\mathcal{G}$ is undirected, connected and time-invariant. \QEDB
\end{assumption}
%
%
%
%
\subsection{Nash Equilibrium Seeking Dynamics}
In order to achieve NES with fixed-time $\mathcal{K}\mathcal{L}_{\mathcal{T}}$ bounds, each player updates its own action using the feedback rule
\begin{equation}\label{additive_dither}
u_i=\hat{u}_i+a_i\hat{\mu}_i,
\end{equation}
where $a_i\in(0,1)$ is a tunable parameter, and $\hat{u}_i$ is updated according to the following dynamics
\begin{equation}\label{fixed_time_dynamics}
\dot{\hat{u}}_i=-kx_{ii}\left(\frac{1}{\left(x_i^\top x_i\right)^\frac{\alpha_{1}}{2}}+\frac{1}{\left(x^\top_i x_i\right)^\frac{\alpha_{2}}{2}}\right),
\end{equation}
where $k>0$ is a tunable gain, and $\alpha_{1},\alpha_{2}$ are defined as
\begin{equation}\label{alphaconstants}
\alpha_{1}:=\frac{q_{1}-2}{q_{1}-1},~~~~\alpha_{2}:=\frac{q_{2}-2}{q_{2}-1},
\end{equation}
where $(q_{1},q_{2})\in\mathbb{R}_{>0}^2$ are tunable parameters that are said to be \emph{admissible} if they satisfy $q_{1}\in(2,\infty)$ and $q_{2}\in(1,2)$. Admissible parameters $(q_{1},q_{2})$ guarantee that $\alpha_{1}$ is always positive and $\alpha_{2}$ is always negative. To implement \eqref{fixed_time_dynamics}, each player is endowed with a vector state $x_i:=[x_{i1},x_{i2},\ldots,x_{i,N}]^\top\in\mathbb{R}^N$, with dynamics $\dot{x}_{ij}$ given by
\begin{equation}\label{fast_consensus}
\varepsilon_{1,i}\dot{x}_{ij}=\sum_{k\in\mathcal{N}_i}\bigg(x_{kj}-x_{ij}\bigg)+b_{ij}\left(\frac{2}{a}J_i(u)\hat{\mu}_i-x_{ij}\right),
\end{equation}
where $\varepsilon_{1,i}:=\varepsilon_1\rho_{1,i}$ is a parameter of each player, $1\gg\varepsilon_{1}>0$ is a tunable parameter, and $b_{ij}=1$ if $i=j$, and $b_{ij}=0$ for all $i\neq j$. For simplicity, in this paper we assume that $\rho_{1,i}=1$ for all $i$. The feedback law \eqref{additive_dither} and the dynamics \eqref{fast_consensus}  depend also on an  auxiliary state $\hat{\mu}_i$, which is generated by each player as the first component of the solution of a linear oscillator with state $\bar{\mu}_i:=[\hat{\mu}_i,\tilde{\mu}_i]^\top\in\mathbb{S}^1$:
\begin{equation}\label{oscillator_dynamics0}
\varepsilon_{2,i}\left[\begin{array}{c}
\dot{\hat{\mu}}_i\\
\dot{\tilde{\mu}}_i
\end{array}\right]=2\pi \mathcal{R}_i\left[\begin{array}{c}
\hat{\mu}_i\\
\tilde{\mu}_i
\end{array}\right],~~\mathcal{R}_i:=\left[\begin{array}{cc}
0 & -\kappa_i\\
\kappa_i & 0
\end{array}\right],
\end{equation}
where $\varepsilon_{2,i}=\varepsilon_{2}\rho_{2,i}$ is a parameter of each player, $\varepsilon_{1}\gg \varepsilon_{2}>0$, $\rho_{2,i}>0$, and $\kappa_i>0$ are tunable parameters. We make the following assumption on the parameters. 
\begin{assumption}\label{frequencies}
Let $\tilde{\kappa}_i:=\kappa_i/\rho_{2,i}$. For each player $i\in\mathcal{V}$, the parameter $\tilde{\kappa}_i$ is a positive rational number, $\tilde{\kappa}_i\neq\tilde{\kappa}_j$, and $\tilde{\kappa}_i\neq2\tilde{\kappa}_j$ for all $i\neq j$. \QEDB
\end{assumption}
%
%
\subsection{FxTNES in Potential Games}
To study the FxTNES dynamics, we start by considering a class of games  termed \emph{potential games} \cite{PotentialGames}, which are characterized by the following assumption.
\begin{assumption}\label{potential_game}
There exists a $\mathcal{C}^2$ radially unbounded function $P:\mathbb{R}^N\to\mathbb{R}$ satisfying $\frac{\partial P(u)}{\partial u_i}=\frac{\partial J_i(u_i,u_{-i})}{\partial u_i}$, for all $u\in\mathbb{R}^N$ and all $i\in\mathcal{V}$; and there exists a unique $u^*\in\mathbb{R}^N$ such that $\frac{\partial J_i(u^*_i,u^*_{-i})}{\partial u_i}=0$, for all $i\in\mathcal{V}$, and $u^*=\arg \min_{u\in\mathbb{R}^N}P(u)$. \QEDB
\end{assumption}

We will also make the following assumption on the mapping $u\mapsto P(u)$.
\begin{assumption}\label{assumption_KL}
There exists $\kappa>0$ such that $P(u)-P(u^*)\leq  \frac{1}{2\kappa}|\nabla P(u)|^2$, for all $u\in\mathbb{R}^N$. \QEDB
\end{assumption}

To characterize the fixed-time convergence properties of the Nash seeking dynamics in potential games, we introduce the constants $\gamma_1=2^{\frac{8-3\alpha_1}{4}}\kappa^{\frac{2-\alpha_1}{2}}$, and $\gamma_2=2^{\frac{8-3\alpha_2}{4}}\kappa^{\frac{2-\alpha_2}{2}}$. Using these constants, we define the following fixed-time:
\begin{equation}\label{fixed_time1}
T_P^*:=\frac{4}{k}\left(\frac{1}{\gamma_1\alpha_1}-\frac{1}{\gamma_2\alpha_2}\right).
\end{equation}
When $(q_{1},q_{2})$ are admissible, the term inside the parenthesis of \eqref{fixed_time1} is positive. Thus, for any desired $T_P^*>0$, equation \eqref{fixed_time1} can be satisfied by simply selecting the gain $k$ of the players as $k=\frac{4}{T_P^*}(\frac{1}{\gamma_1\alpha_1}-\frac{1}{\gamma_2\alpha_2})$. The following theorem is the first result of the paper.
\begin{theorem}\label{theorem1}
Consider the FxTNES dynamics given by \eqref{fixed_time_dynamics}, \eqref{fast_consensus}, and \eqref{oscillator_dynamics0}. Suppose that Assumptions \ref{graph_assumption}, \ref{frequencies}, \ref{potential_game}, and \ref{assumption_KL} hold. Then, for any $T_P^*>0$, there exist admissible parameters $(q_1,q_2)$ and gain $k>0$, such that there exists $\beta\in\mathcal{K}\mathcal{L}_{\mathcal{T}}$ such that the following holds: 
\begin{enumerate}
\item For each pair $\delta>\nu>0$, there exists $\varepsilon_1^*>0$ such that for each $\varepsilon_1\in(0,\varepsilon_1^*)$ there exists $a^*>0$ such that for each $a\in(0,a^*)$ there exists $\varepsilon_2^*>0$ such that for each $\varepsilon_2\in(0,\varepsilon_2^*)$ the FxTNES dynamics with initial conditions $(\hat{u}_i(0),x_i(0),\bar{\mu}_i(0))\in \left(\{u^*\}+\delta\mathbb{B}\right)\times \delta\mathbb{B}\times\mathbb{T}^n$ generate complete solutions, and the vector of actions of the players satisfies the bound 
\begin{equation}\label{KLbound1}
|u(t)-u^*|\leq \beta(|\hat{u}(0)-u^*|,t)+\nu,~~\forall~t\geq0.
\end{equation}
\item $\beta(r,s)=0$ for all $s\geq T^*_{P}$ and all $r>0$; i.e., $\beta$ has the fixed-time convergence property.   \QEDB
\end{enumerate}
\end{theorem}

The dependence of the parameters $(a,\varepsilon_1,\varepsilon_2)$ on the pair $(\delta,\nu)$ implies that Theorem \ref{theorem1} is a semi-global practical asymptotic stability result. However, unlike previous results in the literature, e.g., \cite{Frihauf12a}, the convergence of the actions of the players is characterized by a class $\mathcal{K}\mathcal{L}_{\mathcal{T}}$ function with a fixed-time convergence property, where $T^*_{P}$ can be completely prescribed \emph{a priori} by the system designer, provided a lower bound on $\kappa$ is known. In particular, the gains $k_i$ do not need to be retuned in order to maintain the convergence bound $T_P^*$ as $\delta$ (resp. $\nu$) increases (resp. decreases). 
\begin{remark}
As shown later in the proof of Theorem \ref{theorem1}, the semi-global practical nature of the result and the use of singular perturbation theory (where \eqref{fast_consensus} and \eqref{oscillator_dynamics0} act as fast dynamics) allow us to consider exponentially stable consensus fast dynamics of the form \eqref{fast_consensus}, instead of fixed-time consensus dynamics \cite{fixedtimebook}, for which the application of averaging theory might be intractable in our case. \QEDB
\end{remark}
\begin{remark}
Even though the dynamics \eqref{fixed_time_dynamics} use homogenous gains $k$ and homogenous exponents ($\alpha_1,\alpha_2$), it is possible to obtain similar convergence results using heterogenous exponents and gains. Such results are omitted in this paper due to space limitations.
\end{remark}
\vspace{-0.1cm}
\subsection{FxTNES in Strongly Monotone Games}
\label{sec_strong}
We now consider general non-cooperative games characterized by $\kappa$-strongly monotone pseudo-gradients $G$ obtained from the vector of costs $J(u)=[J_1(u),J_2(u),\ldots,J_N(u)]^\top$. These games, termed \emph{strongly monotone games}, are characterized by the following assumption.
\begin{assumption}\label{assumption_strong_monotonicity}
There exists a unique Nash equilibrium $u^*\in\mathbb{R}^N$, and the pseudo-gradient mapping $u\mapsto G(u)$ satisfies the following inequality $(u_1-u_2)^\top \left(G(u_1)-G(u_2)\right)\geq\kappa |u_1-u_2|^2$, for all $u_1,u_2\in\mathbb{R}^n$, and for $\kappa>0$. \QEDB
\end{assumption}
\begin{remark}
Strongly monotone games have been extensively studied in the literature of variational inequalities \cite{FVIP1}, and distributed control \cite{PersisGrammatico}. \QEDB
\end{remark}
%

For strongly monotone games, we define the constants $\theta_1=2^{-\frac{\alpha_1}{2}}$, $\theta_2=2^{-\frac{\alpha_2}{2}}$, where $\alpha_{1}$ and $\alpha_{2}$ are defined in \eqref{alphaconstants}. Using these constants, we define the fixed-time:
\begin{equation}\label{bound_monotone}
T_S^*:= \frac{4}{k\kappa}\left(\frac{\theta_1}{\alpha_1}-\frac{\theta_2}{\alpha_2}\right).
\end{equation}
Since for any admissible pair $(q_{1},q_{2})$ the term inside the parentheses is positive, for any $T_S^*>0$ we can satisfy equation \eqref{bound_monotone} by setting  $k:= \frac{4}{T_S^*\kappa}\left(\frac{\theta_1}{\alpha_1}-\frac{\theta_2}{\alpha_2}\right)$. The following theorem corresponds to the second main result of the paper.
\begin{theorem}\label{theorem2}
Consider the FxTNES dynamics given by \eqref{fixed_time_dynamics}, \eqref{fast_consensus}, and \eqref{oscillator_dynamics0}. Suppose that Assumptions \ref{graph_assumption}, \ref{frequencies} and \ref{assumption_strong_monotonicity} hold. Then, for any $T_S^*>0$, there exist admissible parameters $(q_1,q_2)$ and gains $k_i>0$, for all $i\in\mathcal{V}$, such that there exists $\beta\in\mathcal{K}\mathcal{L}_{\mathcal{T}}$ such that the following holds: 
\begin{enumerate}
\item For each pair $\delta>\nu>0$ there exists $\varepsilon_1^*>0$ such that for each $\varepsilon_1\in(0,\varepsilon_1^*)$ there exists $a^*>0$ such that for each $a\in(0,a^*)$ there exists $\varepsilon_2^*>0$ such that for each $\varepsilon_2\in(0,\varepsilon_2^*)$ the FxTNES dynamics with initial conditions $(\hat{u}_i(0),x_i(0),\bar{\mu}_i(0))\in \left(\{u^*\}+\delta\mathbb{B}\right)\times \delta\mathbb{B}\times\mathbb{T}^n$ generate solutions with unbounded time domain, and each of these solutions satisfies 
\begin{equation}\label{KLbound2}
|u(t)-u^*|\leq \beta(|\hat{u}(0)-u^*|,t)+\nu,~\forall~t\geq0.
\end{equation}
\item $\beta(r,s)=0$ for all $s\geq T^*_{S}$ and all $r>0$; i.e., $\beta$ has the fixed-time convergence property.    \QEDB
\end{enumerate}
\end{theorem}

\vspace{0.1cm}

\section{Stability and Convergence Analysis}
\label{sec_proofs}
In this section, we present the analysis and proofs of our main results.  Since the FxTNES dynamics are continuous but not Lipschitz continuous, standard averaging and singular perturbation tools for smooth extremum seeking dynamics cannot be used as in \cite{Frihauf12a}. Instead, we will use generalized averaging tools for extremum seeking dynamics that only use continuity of the dynamics, e.g., \cite{WangAveraging,PoTe16,zero_order_poveda_Lina}. 
We organize the proof of the Theorems into multiple steps.

\vspace{0.1cm}
\noindent
\textsl{Step 0: Setting up the Model.} Let us start by writing the FxTNES dynamics as a singularly perturbed system \cite{WangAveraging,zero_order_poveda_Lina}. To do this, define the vectors $\mathbf{x}_{ii}:=[x_{11},x_{22},x_{33},\ldots, x_{NN}]^\top\in\mathbb{R}^N$, as well as the vector-valued mapping $\Psi(\mathbf{x}):=\left[\psi(\mathbf{x}_1),\psi(\mathbf{x}_2),\ldots,\psi(\mathbf{x}_N)\right]^\top\in\mathbb{R}^N$, with mappings $z\to\psi(z)$ defined as
\begin{equation}\label{individual_mappings}
\psi(z):=(z^\top z)^{-\frac{\alpha_{1}}{2}}+(z^\top z)^{-\frac{\alpha_{2}}{2}},
\end{equation}
for all $z\in\mathbb{R}^N$. We also define the vectors
\begin{align*}
\mathbf{x}_i:&=[x_{i1},x_{i2},x_{i3},\ldots,x_{iN}]^\top\in\mathbb{R}^N,\\
\mathbf{x}:&=[\mathbf{x}^\top_1,\mathbf{x}^\top_2,\mathbf{x}^\top_3,\ldots,\mathbf{x}^\top_N]^\top\in\mathbb{R}^{NN}.\\
\mathbf{b}_i:&=[b_{i1},b_{i2},b_{i3},\ldots,b_{iN}]^\top\in\mathbb{R}^N,\\
\mathbf{b}:&=[\mathbf{b}^\top_1,\mathbf{b}^\top_2,\mathbf{b}^\top_3,\ldots,\mathbf{b}^\top_N]^\top\in\mathbb{R}^{NN},\\
\mathbf{\mu}:&=[\hat{\mu}_1,\tilde{\mu}_1,\hat{\mu}_2,\tilde{\mu}_2,\ldots,\hat{\mu}_N,\tilde{\mu}_N]^\top\in\mathbb{R}^{2N},
\end{align*}
and the matrix $B:=\text{diag}(\mathbf{b})$. Let 
\begin{align*}
\mathbf{J}(u,\hat{\mu})&:=[J_1(u)\hat{\mu}_1,J_2(u)\hat{\mu}_2,\ldots,J_N(u)\hat{\mu}_N]^\top\in\mathbb{R}^N,\\
A&:=\text{diag}([a_1,a_2,a_3,\ldots,a_N]^\top)\in\mathbb{R}^{N\times N},
\end{align*}
and define the mapping $E(u,\hat{\mu}):=2(A^{-1}\mathbf{J}(u,\hat{\mu}))\otimes\mathbf{1_N}$. Finally, define a block diagonal matrix $\mathcal{R}_{\tilde{\kappa}}\in\mathbb{R}^{2N\times 2N}$, which is parametrized by a vector of gains $\tilde{\kappa}:=[\tilde{\kappa}_1,\tilde{\kappa}_2,\ldots,\tilde{\kappa}_N]^\top$ (c.f., Assumption \ref{frequencies}), with the $i^{th}$ diagonal block of $\mathcal{R}_{\tilde{\kappa}}$ defined as \eqref{oscillator_dynamics0}.

By using these definitions, the FxTNES dynamics can be written in vectorial form as
\begin{subequations}\label{flow_map}
\begin{align}
\dot{\hat{u}}&=-k\text{diag}(\mathbf{x}_{ii})\Psi(x)\label{learning_dynamics}\\
\varepsilon_1\dot{\mathbf{x}}&=-(L\otimes I_{N\times N}+B)\mathbf{x}+BE(\hat{u}+A\hat{\mu})\label{consensus_dynamics}\\
\varepsilon_2\dot{\mathbf{\mu}}&=2\pi R_{\tilde{\kappa}}\mathbf{\mu}.\label{oscillator_dynamics}
\end{align}
\end{subequations}
For the purpose of algorithmic analysis, we will initially restrict the states $(\hat{u},\mathbf{x},\mu)$ to evolve in the following sets
\begin{equation}\label{flow_set}
(\hat{u},\mathbf{x},\mu)\in \mathbb{R}^N\times M\mathbb{B}\times\mathbb{T}^N,
\end{equation}
where $M>0$ is a constant that can be taken arbitrarily large. Equations \eqref{flow_map} and \eqref{flow_set} describe a constrained singularly perturbed dynamical system, see \cite{WangAveraging,zero_order_poveda_Lina}. Note that, by construction, the set $\mathbb{T}^N$ is forward invariant under the dynamics of $\mu$. The following lemma follows the same ideas of \cite[Sec. III-B]{JournalPovedaKrstic}. 
\begin{lemma}
Suppose the parameters $q_{1}$ and $q_{2}$ are admissible for every player $i\in\mathcal{V}$. Then, the right-hand side of \eqref{flow_map} is continuous in $(\hat{u},\mathbf{x},\mu)$.
\end{lemma}

\vspace{0.1cm}
\noindent
\textsl{Step 1: First Application of Averaging Theory}. We now proceed to analyze system \eqref{flow_map}-\eqref{flow_set} using averaging theory. The following lemma will be instrumental for our results. The proof can be found in \cite[Appendix A]{zero_order_poveda_Lina}.
\begin{lemma}\label{properties_oscillator}
Suppose that Assumption \ref{frequencies} holds. Then, there exists a $T>0$ such that every solution $\mu$ of the oscillator \eqref{oscillator_dynamics} with $\varepsilon_2=1$ satisfies $\frac{1}{\ell T}\int_{0}^{\ell T}\hat{\mu}(t)\hat{\mu}(t)^\top dt=\frac{1}{2}I_n$ and $\frac{1}{\ell T}\int_{0}^{\ell T}\hat{\mu}(t)dt=\textbf{0}$, for all $\ell\in\mathbb{Z}_{\geq1}$, where $\hat{\mu}=[\hat{\mu}_1,\hat{\mu}_2,\ldots,\hat{\mu}_N]^\top$. \QEDB
\end{lemma}

Using Lemma \ref{properties_oscillator}, we will now proceed to average the dynamics \eqref{learning_dynamics}-\eqref{consensus_dynamics} along the solutions of system \eqref{oscillator_dynamics}. In order to do this, and for values of $a_i>0$ sufficiently small, we consider a Taylor expansion of each cost function $J_i(\hat{u}+A\hat{\mu})$ around $u$, leading to $J_i(\hat{u}+A\hat{\mu})=J_i(\hat{u})+\sum_{k=1}^Na_k\hat{\mu}_k\frac{\partial J_i(\hat{u})}{\partial \hat{u}_k}+O_i(a^2)$, where $O_i(a^2)$ represents higher order terms that are bounded on compact sets, and which can be made arbitrarily small by decreasing $a_i$. Using Lemma \ref{properties_oscillator}, it follows that $\frac{1}{\ell T}\int_{0}^{\ell T}\hat{\mu}_i(s)J_i(\hat{u}+A\hat{\mu}(s))ds=\frac{a_i}{2}\frac{\partial J_i(\hat{u})}{\partial u_i}+O_i(a^2)$, for each $i\in\mathcal{V}$. It then follows that $\frac{1}{\ell T}\int_{0}^{\ell T} E\left(\hat{u}+A\hat{\mu}(s),\hat{\mu}(s)\right)ds=G(\hat{u})\otimes\mathbf{1_N}+O(a)$, where $G$ is the pseudo-gradient of the game. Substituting in \eqref{consensus_dynamics}, we obtain the following average dynamics with state $(\hat{u}^a,\mathbf{x}^a)$:
\begin{subequations}\label{ES_dynamics_average}
\begin{align}
\dot{\hat{u}}^a&=-k\text{diag}(\mathbf{x}^a_{ii})\Psi(x^a)\\
\varepsilon_1\dot{\mathbf{x}}^a&=-(L\otimes I_{N\times N}+B)\mathbf{x}^a+BG(\hat{u}^a)\otimes\mathbf{1_N}+O(a),\label{dynamics_linearfilter}
\end{align}
\end{subequations}

\vspace{-0.4cm}
\noindent which evolve in the set
\begin{equation}\label{set_average1}
(\hat{u}^a,\mathbf{x}^a)\in \mathbb{R}^N\times M\mathbb{B}.
\end{equation}
System \eqref{ES_dynamics_average}-\eqref{set_average1} is an $O(a)$-perturbed version of a \emph{nominal average} system with dynamics
\begin{subequations}\label{ES_dynamics_average2}
\begin{align}
\dot{\hat{u}}^a&=-k\text{diag}(\mathbf{x}^a_{ii})\Psi(x^a)\label{learning_a}\\
\varepsilon_1\dot{\mathbf{x}}^a&=-(L\otimes I_{N\times N}+B)\mathbf{x}^a+BG(\hat{u}^a)\otimes\mathbf{1_N}.
\end{align}
\end{subequations}
Therefore, we proceed to analyze the stability and convergence properties of the nominal system \eqref{set_average1}-\eqref{ES_dynamics_average2}, using  robustness results for continuous ODEs to establish stability properties for the perturbed dynamics \eqref{ES_dynamics_average}.

\vspace{0.1cm}
\noindent
\textsl{Step 2: Second Application of Averaging Theory.}
For $\varepsilon_1$ sufficiently small, system  \eqref{set_average1}-\eqref{ES_dynamics_average2} is in singular perturbation form with the dynamics of $\mathbf{x}^a$ acting as fast dynamics, see \cite{WangAveraging,zero_order_poveda_Lina}. To find the boundary layer dynamics, let $\tau=t/\varepsilon_1$,  and consider the system in the $\tau$-time scale:
\begin{subequations}\label{ES_dynamics_average_tau}
\begin{align}
\frac{\partial\hat{u}}{\partial \tau}^a&=-\varepsilon_1k\text{diag}(\mathbf{x}^a_{ii})\Psi(x^a)\\
\frac{\partial \mathbf{x}^a}{\partial \tau}&=-(L\otimes I_{N\times N}+B)\mathbf{x}^a+BG(\hat{u}^a)\otimes\mathbf{1_N}.
\end{align}
\end{subequations}
Setting $\varepsilon_1=0$, we obtain the boundary layer dynamics
\begin{align}\label{ES_dynamics_average_bl}
\frac{\partial \mathbf{x}_{bl}^a}{\partial \tau}=-(L\otimes I_{N\times N}+B)\mathbf{x}^a_{bl}+BG(\hat{u}^a_{bl})\otimes\mathbf{1_N},
\end{align}
where $u^a_{bl}$ is constant. Since $b_{ii}=1$ for all $i\in\mathcal{V}$, the matrix $-(L\otimes I_{N\times N}+B)$ is Hurwitz \cite[Lemma 1.6]{MASBook11}, and since the entries of $B$ are either $0$ or $1$, we obtain that the states $x^a_{bl,ij}$ of system \eqref{ES_dynamics_average_bl} converge exponentially fast to the equilibria $x_{bl,ij}^{a*}=\frac{\partial J_j(\hat{u}^a_{bl})}{\partial \hat{u}^a_{bl,j}}$, for all $i,j\in\mathcal{V}$. In particular, note that $$\mathbf{x}_{bl,i}^{a*\top} \mathbf{x}^{a*}_{bl,i}= G(\hat{u}_{bl}^a)^\top G(\hat{u}_{bl}^a)=|G(\hat{u}_{bl}^a)|^2,$$ for all $i\in\mathcal{V}$. Thus, the singularly perturbed system \eqref{set_average1}-\eqref{ES_dynamics_average2} has a well-defined reduced system, see \cite[Ex. 1]{WangAveraging}, which corresponds to the dynamics \eqref{learning_a} with $x^a$ substituted by the steady-state values $x_{bl,ij}^{a,*}$. Therefore, the reduced system, with state $z\in\mathbb{R}^N$, are given by $\dot{z}=-k\text{diag}(G(z))\Upsilon(z)$, where $\Upsilon$ is defined as $\Upsilon(z):=\left[\psi(G(z)),\psi(G(z)),\ldots,\psi(G(z))\right]^\top\in\mathbb{R}^N$, with individual mappings $z\to\psi(z)$ defined as in \eqref{individual_mappings}. It follows that the reduced dynamics of each agent $i$ are given by
\begin{align}\label{ES_dynamics_average_nominal_reduced1}
\dot{z}_i=-k\frac{\partial J_i(z)}{\partial z_i}\left(\dfrac{1}{|G(z)|^{\alpha_{1}}}+\dfrac{1}{|G(z)|^{\alpha_{2}}}\right).
\end{align}

\vspace{0.1cm}
\noindent
\textsl{Step 3:  Fixed-Time Convergence $\mathcal{K}\mathcal{L}_{\mathcal{T}}$ bounds.} Next, we establish suitable fixed-time $\mathcal{K}\mathcal{L}_{\mathcal{T}}$ convergence bounds for the dynamics \eqref{ES_dynamics_average_nominal_reduced1}. 
\begin{lemma}\label{lemma1}
Suppose that Assumptions \ref{potential_game} and \ref{assumption_KL} hold, and consider the dynamics \eqref{ES_dynamics_average_nominal_reduced1} for all $i\in\mathcal{V}$. Then, there exists $\beta\in\mathcal{K}\mathcal{L}_{\mathcal{T}}$ such that every solution satisfies $|z(t)-u^*|\leq\beta(|z(0)-u^*|,t)$, for all $t\geq0$, and $\beta(r,s)=0$ for all $s>T_P^*$, and all $r>0$, where  $T^*_P$ is given by \eqref{fixed_time1}. \QEDB
\end{lemma}

\textsl{Proof:} Consider the Lyapunov function $V_P(z)=\frac{1}{2}(P(z)-P(u^*))^2$, which under Assumption \ref{potential_game} is positive definite with respect to $u^*$, and also radially unbounded. 
%
Since the game is a potential game, we have that $\nabla P(z)=G(z)$, and the dynamics \eqref{ES_dynamics_average_nominal_reduced1} can be written as
\begin{align*}
\dot{z}=-kG(z)\left(\dfrac{1}{|G(z)|^{\alpha_{1}}}+\dfrac{1}{|G(z)|^{\alpha_{2}}}\right).
\end{align*}
Thus, the time-derivative of $V$ satisfies
\begin{align*}
\dot{V}_P(z)=-k\big(P(z)-P^*\big)\left(\dfrac{G(z)^\top G(z)}{|G(z)|^{\alpha_{1}}}+\dfrac{G(z)^\top G(z)}{|G(z)|^{\alpha_{2}}}\right).
\end{align*}
Defining $\tilde{\alpha}_1:=2-\alpha_1>0$, $\tilde{\alpha}_2:=2-\alpha_2>0$, $c_1:=2^{\frac{2+3\tilde{\alpha_1}}{4}}\kappa^{\frac{\tilde{\alpha}_1}{2}}>0$, $c_2:=2^{\frac{2+3\tilde{\alpha_2}}{4}}\kappa^{\frac{\tilde{\alpha}_2}{2}}>0$, and using Assumption \ref{assumption_KL} we obtain 
\begin{align*}
\dot{V}_P(z)\leq -k \Big(c_1V_P(z)^{\gamma_1}+c_2V_P(z)^{\gamma_2} \Big)<0,~~~\forall~z\neq u^*,
\end{align*}
where  $\gamma_1:=\frac{2+\tilde{\alpha}_1}{4} \in(0,1)$,$\gamma_2:=\frac{2+\tilde{\alpha}_2}{4}>1$. The result follows by \cite[Lemma 1]{Fixed_timeTAC}.  \hfill $\blacksquare$

\begin{lemma}\label{lemma2}
Suppose that Assumption \ref{assumption_strong_monotonicity} holds and that every player implements the dynamics \eqref{ES_dynamics_average_nominal_reduced1}. Then, there exists $\beta\in\mathcal{K}\mathcal{L}_{\mathcal{T}}$ such that every solution satisfies $|z(t)-u^*|\leq\beta(|z(0)-u^*|,t)$, for all $t\geq0$, and $\beta(r,s)=0$ for all $s>T_S^*$, and all $r>0$, where $T^*_S$ is given by \eqref{bound_monotone}. \QEDB
\end{lemma}

\textsl{Proof:} Consider the Lyapunov function $V(z)=\frac{1}{2}|G(z)|^2$, which is positive definite with respect to the Nash equilibrium $u^*$, and also radially unbounded due to the strong monotonicity assumption. The time derivative of $V$ satisfies $\dot{V}(z)=G(z)^\top \mathcal{J}G(z)\dot{z}$ where $\mathcal{J}G$ is the Jacobian matrix of the pseudo-gradient $G$. Since $G$ is strongly $\kappa$-monotone, $\eta^\top \mathcal{J}G(z)\eta\geq \kappa |\eta|^2$,~~for all $\eta\in\mathbb{R}^n$, and for all $z\in\mathbb{R}^n$ \cite[Prop. 2.3.2]{FVIP1}. Therefore,  we obtain
\begin{align*}
\dot{V}(z)&\leq-k\kappa  \left( c_1V(z)^{\gamma_1}+c_2V(z)^{\gamma_2}  \right)<0,~~~\forall~z\neq u^*,
\end{align*}
where $c_1=2^{\frac{\tilde{\alpha}_1}{2}}$, $c_2=2^{\frac{\tilde{\alpha}_2}{2}}$, $\gamma_1=\frac{\tilde{\alpha}_1}{2}$,$\gamma_2=\frac{\tilde{\alpha}_2}{2}$. The result follows by \cite[Lemma 1]{Fixed_timeTAC}. \hfill $\blacksquare$

\vspace{0.1cm}
\noindent
\textsl{Step 4: $\beta$-SGPAS for Second Singularly Perturbed System.} Having established UGAS of the average dynamics with $\mathcal{K}\mathcal{L}_{{\mathcal{T}}}$ bounds $\beta$, we proceed to apply averaging results in order to establish suitable stability properties for the singularly perturbed dynamics \eqref{set_average1}-\eqref{ES_dynamics_average2}. In particular, by \cite[Thm. 2]{WangAveraging} the dynamical system  \eqref{set_average1}-\eqref{ES_dynamics_average2} renders the compact set  $\mathcal{A}:=\{u^*\}\times M\mathbb{B}$, $\beta$-SGPAS as $\varepsilon_1\to 0^+$ with the $\mathcal{K}\mathcal{L}_{\mathcal{T}}$ bound $\beta$ obtained in Lemmas \ref{lemma1} or \ref{lemma2}.  Moreover, by the definition of solutions we have that $|\mathbf{x}^a(t)|_{M\mathbb{B}}=0$ for all $t\in\text{dom}(u^a,\mathbf{x}^a)$, which implies that $|\zeta^a(t)|_{\mathcal{A}}:=|u^a(t)-u^*|$, where $\zeta^a:=[u^{a\top},\mathbf{x}^{a\top}]^\top$. Thus, for each $\delta>\nu>0$, there exists $\varepsilon_1^*>0$ such that for all $\varepsilon_1\in(0,\varepsilon^*_1)$ every solution of the dynamical system  \eqref{set_average1}-\eqref{ES_dynamics_average2} with $\zeta^a(0)\in (\{u^*\}+\delta\mathbb{B})\times \delta\mathbb{B}$ satisfies the bound $|\zeta^a(t)|_{\mathcal{A}}\leq \beta(|\zeta^a(0)|_{\mathcal{A}},t)+\nu$, for all $t\in\text{dom}(\zeta^a)$, where $\beta$ comes from Lemmas \ref{lemma1} and \ref{lemma2}. Finally, by \cite[Prop. A.1]{zero_order_poveda_Lina}, the $O(a)$-perturbed system \eqref{ES_dynamics_average}-\eqref{set_average1} renders the same compact set $\mathcal{A}$ $\beta$-SGPAS as $(a,\varepsilon_1)\to 0^+$ with $\mathcal{K}\mathcal{L}_{\mathcal{T}}$ bound $\beta$.

\vspace{0.1cm}
\noindent
\textsl{Step 5: $\beta$-SGPAS for Original Dynamics.} We now analyze the FxTNES dynamics \eqref{flow_map} based on the properties of the average dynamics \eqref{ES_dynamics_average} studied in Step 4. Since the oscillator \eqref{oscillator_dynamics} renders the set $\mathbb{T}^N$ UGAS, it follows by \cite[Thm. 7]{zero_order_poveda_Lina} that the FxTNES dynamics render the compact set $\mathcal{A}\times\mathbb{T}^n$ $\beta$-SGPAS as $(\varepsilon_2,a,\varepsilon_1)\to 0^+$ with the same $\mathcal{K}\mathcal{L}_{\mathcal{T}}$ bound $\beta$ of Step 4. In particular, this establishes that for each $k>0$, for each tuple of admissible parameters $(q_1,q_2)$, and each pair $\delta>\nu>0$, there exists $\varepsilon_1^*>0$ such that for each $\varepsilon_1\in(0,\varepsilon_1^*)$ there exists $a^*\in(0,\nu/2)$ such that for each $a\in(0,a^*)$ there exists $\varepsilon_2^*>0$ such that for each $\varepsilon_2\in(0,\varepsilon_2^*)$ each solution of the FxTNES dynamics satisfies the bound $|\mathbf{y}(t)|_{\mathcal{A}}\leq \beta(|\mathbf{y}(0)|_{\mathcal{A}},t)+\frac{\nu}{2}$, for all $t\in\text{dom}(\mathbf{y})$, where $\mathbf{y}:=(\hat{u}^\top,\mathbf{x}^\top,\mu^\top)^\top$. Given that by definition of solutions we have that $|\mathbf{x},\mu|_{M\mathbb{B}\times\mathbb{T}^N}=0$ and therefore $|\mathbf{y}|_{\mathcal{A}}=|\hat{u}-u^*|$, it follows that the actions of the players satisfy the bound $|\hat{u}(t)-u^*|\leq \beta(|\hat{u}(0)-u^*|,t)+\frac{\nu}{2}$, for all $t\in\text{dom}(\mathbf{y})$. Let $\bar{a}=\max_{i\in\mathcal{V}}a_i$; using the previous $\mathcal{K}\mathcal{L}_{\mathcal{T}}$ bound, $\bar{a}\in(0,\min\{\nu/2,a^*\})$, the triangle inequality, and the fact that $u=\hat{u}+A\hat{\mu}$,  we obtain inequalities \eqref{KLbound1} and \eqref{KLbound2}. Finally, completeness of solutions of the original FxTNES dynamics \eqref{flow_map} without the restriction $M\mathbb{B}$ follows by the linearity of the dynamics \eqref{dynamics_linearfilter}. In particular, for a bounded input $s$, the solutions of a system of the form $\dot{\mathbf{x}}^a=\frac{1}{\varepsilon_1}\left(F\mathbf{x}^a+Us(t)\right)$ satisfy $|\mathbf{x}^a(t)|\leq c_1\exp(-\frac{\lambda}{\varepsilon_1}t)|\mathbf{x}^a(0)|+c_2\|U\|\sup_{t\geq0} |s(t)|$, for all $t\geq0$, some $c_1,c_2>0$, and all $\varepsilon_1>0$. Since in our case $s=G(\hat{u})$, continuity of $G$ and uniform boundedness of $\hat{u}$ imply boundedness of $G(\hat{u})$, and for $\delta>0$ such that $|\hat{u}(0)-u^*|\leq \delta$ and $|\mathbf{x}^a(0)|\leq\delta$, there exists $M>1$ sufficiently large such that $|\mathbf{x}^a(t)|<0.5M$ for all $t\geq0$. The fact that $|\mathbf{x}(t)|<M$ for all $t\geq0$ follows by $\epsilon$-closeness between $\mathbf{x}^a$ and $\mathbf{x}$, and the exponential bound on $|\mathbf{x}^a(t)|$. \null\hfill\null $\blacksquare$

\section{Numerical Results}
\label{sec_numerical}
\begin{figure}[t!]
 \centering
  \includegraphics[width=0.49\textwidth]{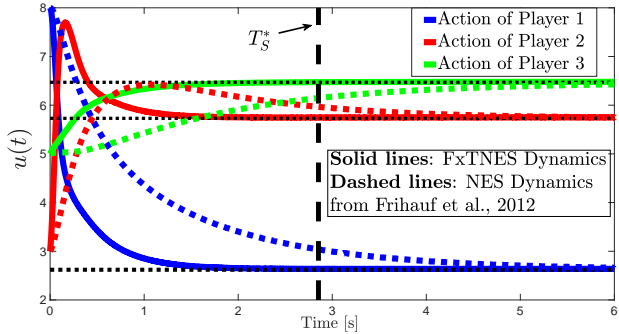}
 \caption{Evolution in time of the actions of the players under the FxTNES dynamics (solid lines) and under the gradient-based NES dynamics considered in \cite{Frihauf12a} (dotted line).}
 \label{figure11}
 \vspace{-0.4cm}
 \end{figure}

\begin{figure}[t!]
 \centering
  \includegraphics[width=0.48\textwidth]{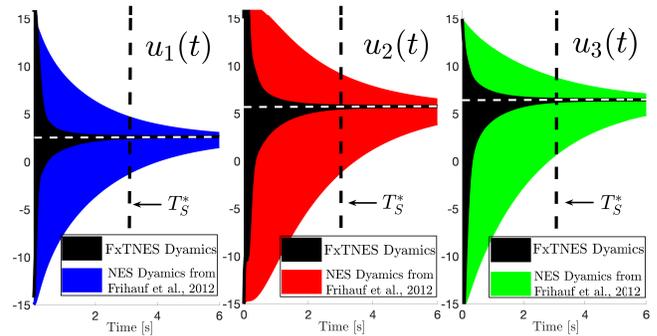}
 \caption{Comparison between the reachable sets in the time interval $[0,6]$ of the NES dynamics of \cite{Frihauf12a} and the FxTNES, with $\hat{u}(0)\in[-15,15]\times[-15,15]\times [-15,15]$.}
 \label{figure22}
 \vspace{-0.4cm}
 \end{figure}

We consider a non-cooperative game with three players having quadratic cost functions of the form $J_i(u)=u^\top Q_i u+b^\top_i u+c_i$. The matrices $Q_i$, the vectors $b_i$, and the constant $c_i$ are selected as $Q_1=[-6, 3, -1;3, 2, 1;-1, 1, 2]$,  $Q_2=
 [3,6,1;6,-9,4;1, 4, 3]$, $Q_3=\left[2,-3,-0.5;-3,-1,1;-0.5, 1,-3\right]$,~$c_i=0$, and $b_1=[10,5,15]^\top,~b_2=[15,20,25]^\top,~b_3=[20,10,30]^\top$. This game is strongly monotone with coefficient $\kappa=4.35$. It also has a unique Nash equilibrium given by $u^*=[2.62, 5.73,6.47]^\top$. In order to achieve model-free NES, we implement the FxTNES dynamics as well as the model-free gradient-based NES dynamics considered in \cite{Frihauf12a}. Both algorithms used the same parameters $k=1$ and $a_i=0.1$ for all $i\in\mathcal{V}$, as well as the constants $\varepsilon_2=1\times 10^{-4}$, and $\varepsilon_1=5\times 10^{-2}$.  
For the FxTNES dynamics, we used the parameters $q_1=3$ and $q_2=1.5$. Figure \ref{figure11} shows the evolution in time of the actions of the players under the FxTNES dynamics (solid lines), and also under the gradient-based NES dynamics of \cite{Frihauf12a} (dotted lines). It can be seen that the trajectories of the FxTNES dynamics converge to the Nash equilibrium before the prescribed time $T_S^*$. To further illustrate the effect of the initial conditions on both algorithms, Figure \ref{figure22} compares the reachable set of both dynamics using the same parameters. In this plot, we can observe a substantial improvement on the transient performance of the trajectories generated by the FxTNES dynamics in comparison to the trajectories generated by the traditional ES gradient-based dynamics .
%
%
\section{Conclusions}
\label{section:conclusions}
We have introduced novel model-free fixed-time Nash equilibrium seeking dynamics for non-cooperative games. In these dynamics, each player needs to evaluate only its own cost function, and to share state information with neighboring players characterized by a communication graph. The convergence of the player's actions to a neighborhood of the Nash equilibrium of the game is dominated by a class $\mathcal{K}\mathcal{L}_{\mathcal{T}}$ function with the ``fixed-time convergence'' property, where the fixed-time can be prescribed a priori by the system designer using minimal information of the game. Numerical examples have illustrated the advantages of the proposed approach compared to traditional gradient descent-based learning rules that have an average system with only asymptotic or exponential convergence properties. Future research directions will focus on dynamic games and settings with time-varying communication topologies.

\bibliographystyle{IEEEtran}
\bibliography{Biblio.bib}

\end{document}